\documentclass{esannV2}
\usepackage[dvips]{graphicx}
\usepackage[latin1]{inputenc}
\usepackage{amssymb,amsmath,array}

\newtheorem{theorem}{Theorem}[section]
\newtheorem{lemma}[theorem]{Lemma}
\begin{document}
%style file for ESANN manuscripts
\title{Consistent estimation of the architecture of multilayer perceptrons}

%***********************************************************************
% AUTHORS INFORMATION AREA
%***********************************************************************
\author{Joseph Rynkiewicz $^1$
%
% Optional short acknowledgment: remove next line if non-needed
%\thanks{This is an optional funding source acknowledgement.}
%
% DO NOT MODIFY THE FOLLOWING '\vspace' ARGUMENT
\vspace{.3cm}\\
%
% Addresses and institutions (remove "1- " in case of a single institution)
Université Paris I - SAMOS/MATISSE\\
90 rue de tolbiac, Paris - France\\ 
}
%***********************************************************************
% END OF AUTHORS INFORMATION AREA
%***********************************************************************

\maketitle

\begin{abstract}
We consider regression models involving multilayer perceptrons (MLP) with one hidden layer and a Gaussian noise. The estimation of the parameters of the MLP can be done by maximizing the likelihood of the model. In this framework,  it is  difficult to determine the true number of hidden units using an information criterion, like the Bayesian information criteria (BIC),  because the information matrix of Fisher is not invertible if the number of hidden units is overestimated. Indeed, the classical theoretical justification of information criteria relies entirely on the invertibility of this matrix. However, using recent methodology introduced to deal with models with a loss of identifiability, we prove that suitable information criterion leads to consistent estimation of the true number of hidden units.
\end{abstract}

\section{Introduction}
Feed-forward neural networks are well known and popular tools to deal with non-linear statistic models.
We can describe MLP regression model as a parametric family of probability density functions. If the noise of the regression model is Gaussian then it is well known (see Watanabe and Fukumizu \cite{Watanabe}) that the maximum likelihood estimator is equal to the least-squares estimator. Therefore, it is natural to consider Gaussian likelihood when we consider feed-forward neural networks from the statistical viewpoint. 
H. White \cite{White} reviews learning in MLP in detail from the statistical viewpoint. However he left pending a important question: The asymptotic behavior of the estimator when an MLP in use has redundant hidden units and the Fisher information matrix is singular. Fukumizu \cite{Fukumizu} gives a response is the case of unbounded parameters. In this case, the maximum likelihood estimator (MLE) behavior is very different from the classical case since the likelihood ratio statistic can have an order lower bounded by $O(\log(n))$ with $n$ the number of observations. This result indicates that the BIC criterion (see Schwarz \cite{Schwarz}) is no more consistent. 

However, it is also a natural assumption to consider that the parameters are bounded. Indeed, computer calculations assume always this boundedness. Moreover it is a safe practice to bound the parameters in order to avoid numerical problems. 

The main result of this paper is to show that if we assume that the parameters are in a suitable compact set (i.e. bounded and closed) the likelihood ratio is tight, so the BIC is convergent.  To obtain this result we use recent techniques introduced by Liu and Shao \cite{Shao} and Gassiat \cite{Gassiat}. These techniques consist in finding a parameterization separating the identifiable part of the parameter vector and the unidentifiable part, then we can obtain an asymptotic development of the likelihood of the model which allows us to show that a set of generalized score functions is a Donsker class. Finally, using a theorem of Gassiat \cite{Gassiat}, we conclude that suitable information criteria like the BIC are consistent because the likelihood ratio statistic is tight. 
\section{The model}
\subsection{Unidentifiability of the true regression function}
Le  $x=(x_1,\cdots,x_d)^T\in{\mathbb R}^d$ be the vector of inputs. The MLP function with $k$ hidden units can be written : 

\[
F_\theta(x)=\beta+\sum_{i=1}^k a_i\phi\left(b_i+w_i^Tx\right)\\
\]
with $\theta=\left(\beta,a_1,\cdots,a_k,b_1,\cdots,b_k,w_{11},\cdots,w_{1d},\cdots,w_{kd}\right)\subset {\mathbb R}^{2k+1+k\times d}$ the parameter vector of the model, and $w_i:=\left(w_{i1},\cdots,w_{id}\right)^T$.
The transfer function $\phi$ will be assumed bounded and three times derivable. We assume also that the first, second and third derivatives of the function $\phi$:  $\phi^{'}$, $\phi^{''}$ and $\phi^{'''}$ are bounded. We consider that the data $\left(X_t,Y_t\right)_{t\in {\mathbb N}^*}$, are random variables verifying the equation: 
\begin{equation}\label{model}
Y_t=F_{\theta^0}(X_t)+\varepsilon_t
\end{equation}
Where $\left(\varepsilon_t\right)_{t\in {\mathbb N}^*}$ is a sequence of independent and identically distributed (i.i.d.) ${\cal N}(0,\sigma^2)$ variables. Note that it is assumed that the true model (\ref{model}) belongs to the considered set of parameter $\Theta$. Define the true number of hidden units as the smallest integer $k^0$ such that it exists \\
$\theta^0=\left(\beta^0,a^0_1,\cdots,a^0_{k^0},b^0_1,\cdots,b^0_{k^0},w^0_{11},\cdots,w^0_{1d},\cdots,w^0_{{k^0}d}\right)\subset {\mathbb R}^{2{k^0}+1+{k^0}\times d}$, with $F_{\theta^0}$ equal to the true regression function of model (\ref{model}). If we overestimate the true number of hidden units then the true parameter will be unidentifiable that is to say it will belong to an union of finitely many submanifolds of $\Theta$ and the dimension of at least one of the manifolds is larger than zero. For example, suppose we have a multilayer perceptron with two hidden units and the true function $F_{\theta^0}(x)$ is given by a perceptron with only one hidden unit, say $F_{\theta^0}(x)=a_1\phi({w^0_1}^Tx)$. then any parameter of the set :
\[
\begin{array}{l}
\left\{\theta\left|w_{1}={w^0_1}, a_1=a^0_1,\beta=b_1=a_2=0\right. \right\}\cup\\
\left\{\theta\left|w_{1}=w_{2}={w^0_1}, a_1+a_2=a^0_1,\beta=b_1=b_2=0\right. \right\}
\end{array}
\]
realizes the function $F_{\theta^0}(x)$. In this framework, the likelihood ratio statistic does not follow the usual chi-square asymptotics, which requires uniqueness of the true parameter in the regularity conditions. An other difficulty appears if it exists a $w_i$ equal to zero, because the function $\phi(b_i+w_i^Tx)$ will be then constant as $\beta$. In order to avoid this source of unidentifiability we will constraint the set of parameters  $\Theta$ to verify for an $\eta>0$, and for all $w_i\in \Theta$:   $\Vert w_i \Vert\geq \eta$.
\subsection{Likelihood of the model}
Let us consider the sample  $Z_i=(X_i,Y_i)$ where $X_i$ and $Y_i$ follow the probability law induced by the model (\ref{model}). We assume that the law of $X_i$ is $q(x)\lambda_d(x)$ with $\lambda_d$ the Lebesgue measure on ${\mathbb R}^d$ and the density function $q(x)$ which is strictly positive for all $x\in{\mathbb R}^d$. The likelihood of the observation $z:=(x,y)$ for a parameter vector $\theta$ will be written: 
\[
f_{\theta}(z)=\frac{1}{\sqrt{2\pi\sigma^2}}e^{-\frac{1}{2\sigma^2}\left(y-F_\theta(x)\right)^2}q(x)
\]
For sake of simplicity and concision we will assume that $\sigma^2$ is known and fixed, but it is not hard to relax this assumption.  
We assume also that it is known that the true number of hidden units is smaller than $M$. $M$ can be very large (for example 1000000), so this assumption is not restrictive in practice. Let \( \Theta:=\cup_{1\leq k\leq M} \Theta_k\) be the set of parameter with an $\eta>0$ such that for all $k$:
\[
\begin{array}{l}
 \Theta_k:=\\
\left\{\theta=\left(\beta,a_1,\cdots,a_k,b_1,\cdots,b_k,w_{11},\cdots,w_{1d},\cdots,w_{kd}\right),\ \forall 1\leq i\leq k, \Vert w_i \Vert\geq \eta\right\}
\end{array}
\]
The set $\Theta$ will be a compact as a finite union of compact sets. We note $k^0$ the true number of hidden units or equivalently the minimal number of hidden units such  that $F_{\theta^0}\in \Theta_{k^0}$ realizes the true regression function. The function $f(z):=f_{\theta^0}(z)$ will be then the true density of the observation. 

\section{Identification of the architecture of the MLP}
Let $l_n(\theta):=\sum_{i=1}^n\log(f_{\theta}(z_i))$ be the log-likelihood of the model, note that this function is known up to the constant $\sum_{i=1}^n\log(x_i)$, independent of the parameter $\theta$. We define $\hat k$, the estimator of maximum of penalized likelihood, as the number of hidden unit maximizing:
\begin{equation}\label{bic}
T_n(k):=\max\{l_n(\theta) : \theta\in\Theta_k\}-p_n(k)
\end{equation}
where $p_n(k)$ is a term which penalizes the log-likelihood in function of the number of hidden units of the model. 
In the sequel, we will assume the following properties:
\begin{description}
\item{H-1 : } The MLP functions are identifiable in the weak following sense: 
\[
\begin{array}{l}
\forall x,\ \beta'+\sum_{i=1}^{k'} a'_i\phi\left(b'_i+{w'_i}^Tx\right)=\beta+\sum_{i=1}^k a_i\phi\left(b_i+w_i^Tx\right)\Leftrightarrow\\
 \delta_{\beta'}+\sum_{i=1}^{k'}a'_i\delta_{(b'_i,w'_i)}=\delta_{\beta}+\sum_{i=1}^{k}a_i\delta_{(b_i,w_i)}
\end{array}
\] where $\delta$ is the Dirac measure, i.e. $\delta_{\theta}(x)=1$ if $x=\theta$ and $\delta_{\theta}(x)=0$ if $x\neq\theta$. 
\item{H-2 :} $E(\left|X\right|^6)<\infty$.
\item{H-3 :} The functions of the set
\[
\begin{array}{l}
\left(\left(x_kx_l\phi^{''}(b^0_i+{w^0_i}^Tx)\right)_{1\leq l \leq k\leq d,\ 1\leq i\leq k^0},\phi^{''}(b^0_i+{w^0_i}^Tx)_{1\leq i\leq k^0}, \right.\\
\left. \left(x_k\phi^{'}(b^0_i+{w^0_i}^Tx)\right)_{1\leq k\leq d,\ 1\leq i\leq k^0}, \left(\phi^{'}(b^0_i+{w^0_i}^Tx)\right)_{1\leq i\leq k^0}\right)
\end{array}
\]
are linearly independents in the Hilbert space $L^2(q\lambda_{d})$. 
\item{H-4 :} $p_n(.)$ is increasing,  $p_n(k_1)-p_n(k_2)\rightarrow \infty$ for all $k_1>k_2$ and $\lim_{n\rightarrow\infty}\frac{p_n(k)}{n}=0$. Note that such conditions are verified by BIC-like criterion. 
\end{description}
We get then the following result:
 
\begin{theorem}
If the assumptions H-1, H-2, H-3 and H-4 are true then  $\hat k\stackrel{P}{\rightarrow}k^0$. 
\end{theorem} 
\paragraph*{Remark}
Sussmann \cite{Sussmann} has shown that, if the transfer functions $\phi$ are sigmoids and if the parameters are  $b_i$ positive (in order to avoid, a symmetry on the signs of $(b_i,w_i)$ and $a_i$), then the assumption (H-1) is verified. Moreover, following a reasoning similar to  Fukimizu \cite{Fukumizu2}, we can show that the sigmoid functions verify the assumption (H-3). So, this result can be applied to the one hidden layer MLP model with sigmoidal transfert functions.
\paragraph*{Sketch of the proof}
Consider the functions:
\[
s_{\theta}(z):=\frac{\frac{f_{\theta}}{f}(z)-1}{\Vert\frac{f_{\theta}}{f}-1\Vert_2} 
\mbox{ where } \Vert .\Vert_2 \mbox{ is the norm } L^2\left(f\lambda_{d+1}\right)
\]
In order to prove the theorem, we have only to show that the set ${\mathbb S}:=\{s_{\theta},\ \theta\in \Theta\}$ is a Donsker class (cf van der Vaart \cite{Vaart}). Roughly speaking, a Donsker class is a set of functions for which the empirical distribution (with i.i.d. variables) verify a uniform central limit theorem, with limit distribution a Gaussian process. Then, the results will follow from the theorem 2.1 of Gassiat \cite{Gassiat}.
Firstly, we will get an asymptotic development of the likelihood ratio when the model is overparametrized. The Donsker property will follow from this development.  
\subsection{Reparameterization of the model}
We will reparameterize the model using the same method as in Liu et Shao \cite{Shao} for the mixing models. 
If $\frac{f_\theta}{f}-1=0$, we have $\beta=\beta^0$ and a vector $t=(t_i)_{1\leq i\leq k^0}$ exists such that $0=t_0<t_1<\cdots<t_{k^0}\leq k$ and up to a permutation:
\(b_{t_{i-1}+1}=\cdots=b_{t_i}=b^0_i\), \(w_{t_{i-1}+1}=\cdots=w_{t_i}=w^0_i\), \(\sum_{j=t_{i-1}+1}^{t_i}a_j=a_i^0\) and \(a_j=0\) for $t_{k^0}+1\leq j\leq k$. Let be $s_i=\sum_{j=t_{i-1}+1}^{t_i}a_j-a_i^0$ and $q_j=\frac{a_j}{\sum_{t_{i-1}+1}^{t_i}a_j}$, we get then the reparameterization $\theta=\left(\Phi_t,\psi_t\right)$ with \(\Phi_t=\left(\beta,(b_j)_{j=1}^{t_{k^0}},(w_j)_{j=1}^{t_{k^0}},(s_i)_{i=1}^{k^0},(a_j)_{j=t_{k^0}+1}^{k}\right)\), \(\psi_t=\left((q_j)_{j=1}^{t_{k^0}},(b_j)_{t_{k^0}+1}^{k},(w_j)_{t_{k^0}+1}^{k}\right)\). With this parameterization, for a fixed $t$,  $\Phi_t$ is an identifiable parameter and all the non-identifiability of the model will be in $\psi_t$. Then, $F_{(\Phi_t^0,\psi_t)}$ will be equal to $F_{\theta^0}$ if and only if

\[
\begin{array}{cccccc}
\Phi^0_t=(\beta^0,&\underbrace{b_1^0,\cdots,b_1^0}&,\cdots,&\underbrace{b_{k^0}^0,\cdots,b_{k^0}^0},&\underbrace{w_1^0,\cdots,w_1^0}&,\cdots,\\
 &t_1& &t_{k^0}-t_{k^0-1}&t_1& \\
&\underbrace{w_{k^0}^0,\cdots,w_{k^0}^0},&\underbrace{0,\cdots,0},&\underbrace{0,\cdots,0})&&\\
&t_{k^0}-t_{k^0-1}&k^0&k-t_{k^0}&&
\end{array}
\]
and now, we have:
\begin{equation}\label{vraisemblance}
\begin{array}{l}
\frac{f_\theta}{f}(z)=\frac{1}{exp\left(-\frac{1}{2\sigma^2}\left(y-\left(\beta^0+\sum_{i=1}^{k^0} a^0_i\phi(b^0_i+{w^0_i}^Tx)\right)\right)^2\right)}\times\\
exp\left(-\frac{1}{2\sigma^2}\left(y-\left(\beta+\sum_{i=1}^{k^0}(s_i+a^0_i)\sum_{j=t_{i-1}+1}^{t_i}q_j\phi(b_j+w_j^Tx)\right.\right.\right.\\
\left.\left.\left.+\sum_{j=t_{k^0}+1}^ka_j\phi(b_j+w_j^Tx)\right)\right)^2\right)
\end{array}
\end{equation}
We get then an approximation of the likelihood ratio  by derivating the expression (\ref{vraisemblance}) with respect to each component of the parameter vector $\phi_t$ and thanks the assumptions H-1, H-2, and H-3.
\begin{lemma}
Let us write $D(\Phi_t,\psi_t):=\Vert\frac{f_{(\Phi_t,\psi_t)}}{f}-1\Vert_2$ and \\
$e(z):=\frac{1}{\sigma^2}\left(y-\left(\beta^0+\sum_{i=1}^{k^0} a^0_i\phi(b^0_i+{w^0_i}^Tx)\right)\right)$ we get the following approximation :
\[
\frac{f_\theta}{f}(z)=1+(\Phi_t-\Phi^0_t)^Tf^{'}_{(\Phi^0_t,\psi_t)}(z)+0.5(\Phi_t-\Phi^0_t)^Tf^{''}_{(\Phi^0_t,\psi_t)}(z)(\Phi_t-\Phi^0_t)+o(D(\Phi_t,\psi_t))
\]
with
\[
\begin{array}{l}
(\Phi_t-\Phi^0_t)^Tf^{'}_{(\Phi^0_t,\psi_t)}(z)=
\left(\beta-\beta^0+\sum_{i=1}^{k^0}s_i\phi(b^0_i+{w^0_i}^Tx)\right.\\
\left. +\sum_{i=1}^{k^0}\sum_{j=t_{i-1}+1}^{t_i}q_j\left(b_j-b^0_i\right)a^0_i\phi^{'}(b^0_i+{w^0_i}^Tx)+\right.\\
\left.\sum_{i=1}^{k^0}\sum_{j=t_{i-1}+1}^{t_i}q_j\left(w_{j}-w^0_{i}\right)^Txa^0_i\phi^{'}(b^0_i+{w^0_i}^Tx)\right.\\
\left.+\sum_{j=t_{k^0}+1}^ka_j\phi(b_j+w_j^Tx)\right)e(z)
\end{array}
\]
and
\[
\begin{array}{l}
(\Phi_t-\Phi^0_t)^Tf^{''}_{(\Phi^0_t,\psi_t)}(z)(\Phi_t-\Phi^0_t)=\\
\left(1-\frac{1}{e^2(z)}\right)\left((\Phi_t-\Phi^0_t)^Tf^{'}_{(\Phi^0_t,\psi_t)}(z){f^{'}_{(\Phi^0_t,\psi_t)}}^T(z)(\Phi_t-\Phi^0_t)\right)\\
+e(z)\times\left(\sum_{i=1}^{k^0}\sum_{j=t_{i-1}+1}^{t_i}q_j(b_j-b^0_i)^2a^0_i\phi^{''}(b^0_i+{w^0_i}^Tx)+\right.\\
\left.\sum_{i=1}^{k^0}\sum_{j=t_{i-1}+1}^{t_i}q_j(w_{j}-w^0_{i})^Txx^T(w_{j}-w^0_{i})a^0_i\phi^{''}(b^0_i+{w^0_i}^Tx)\right.\\
\left. +\sum_{i=1}^{k^0}\sum_{j=t_{i-1}+1}^{t_i}(q_jb_j-b^0_i)s_{i}\phi^{'}(b^0_i+{w^0_i}^Tx)\right.\\
\left.+\sum_{i=1}^{k^0}\sum_{j=t_{i-1}+1}^{t_i}(q_jw_{j}-w^0_{i})^Txs_{i}\phi^{'}(b^0_i+{w^0_i}^Tx) \right)
\end{array}
\]
\end{lemma}
Now, it is easy to show that the minimum number  $N(\epsilon)$ of $\epsilon$-brackets (cf van der Vaart \cite{Vaart}) needed to cover  $\left\{S_\theta, \theta\in\Theta_k\right\}$ is of order $O\left(\frac{1}{\epsilon}^{3k^0}\right)$. This proves that $\mathbb S$ is a Donsker class $\blacksquare$
\section{Conclusion}
Penalized likelihood criterium is used since many years for MLP models. However, even if such technique seems to work in practice, there was no theoritical justification for its use. Indeed, the classical asymptotic theory fails when the Fisher information matrix is singular. For the first time, we give sufficient conditions insuring the success of such selection procedure as the number of observation tends to infinite. This result reinforces the use of classical statistical criterium like BIC in order to fit the architecture of MLP models. 

% ****************************************************************************
% END OF BIBLIOGRAPHY AREA
% ****************************************************************************

%
\end{document}